\documentclass[letterpaper,11pt]{article}
\usepackage{amsmath, amsthm, amssymb}    
\usepackage{bridges}			
\usepackage{graphicx}			
\usepackage[colorlinks=true, urlcolor=blue, citecolor=black, linkcolor=black]{hyperref}  
\usepackage{subcaption}			

\title{Dilute Hitomezashi on the Isometric Grid}

\author{Katherine A. Seaton
\vspace{10pt}\\
La Trobe University, Australia; k.seaton@latrobe.edu.au} 

\date{}					

\begin{document}

\maketitle

\thispagestyle{empty}

\begin{abstract}

This paper describes preliminary investigation of hitomezashi stitching designs created on the isometric grid. An imposed constraint is that only every second line of stitching in each of the possible three directions is present. Each vertex visited by the stitching has degree two. Motifs and wallpaper patterns with three- and six-fold rotational symmetry are illustrated, in particular designs featuring the von Koch snowflake iterates.  

\end{abstract}

\section*{Embroidery on the Isometric Grid}
Hitomezashi is a traditional form of embroidery dating from Edo period Japan wherein perpendicular lines of running stitch with common vertices interact to form motifs and overall designs on the square grid.  Uniform length stitches are placed alternately on the face and the reverse of the fabric being worked upon, with the result that a second design (the dual) forms on what in other forms of stitching would be called the `wrong' side. Since hitomezashi appeared at Bridges in 2020 \cite{HayesSeat}, a growing number of mathematicians and mathematical fibre artists have explored its potential.

\begin{figure}[h!tbp]
	\centering
	\includegraphics[width=5in]{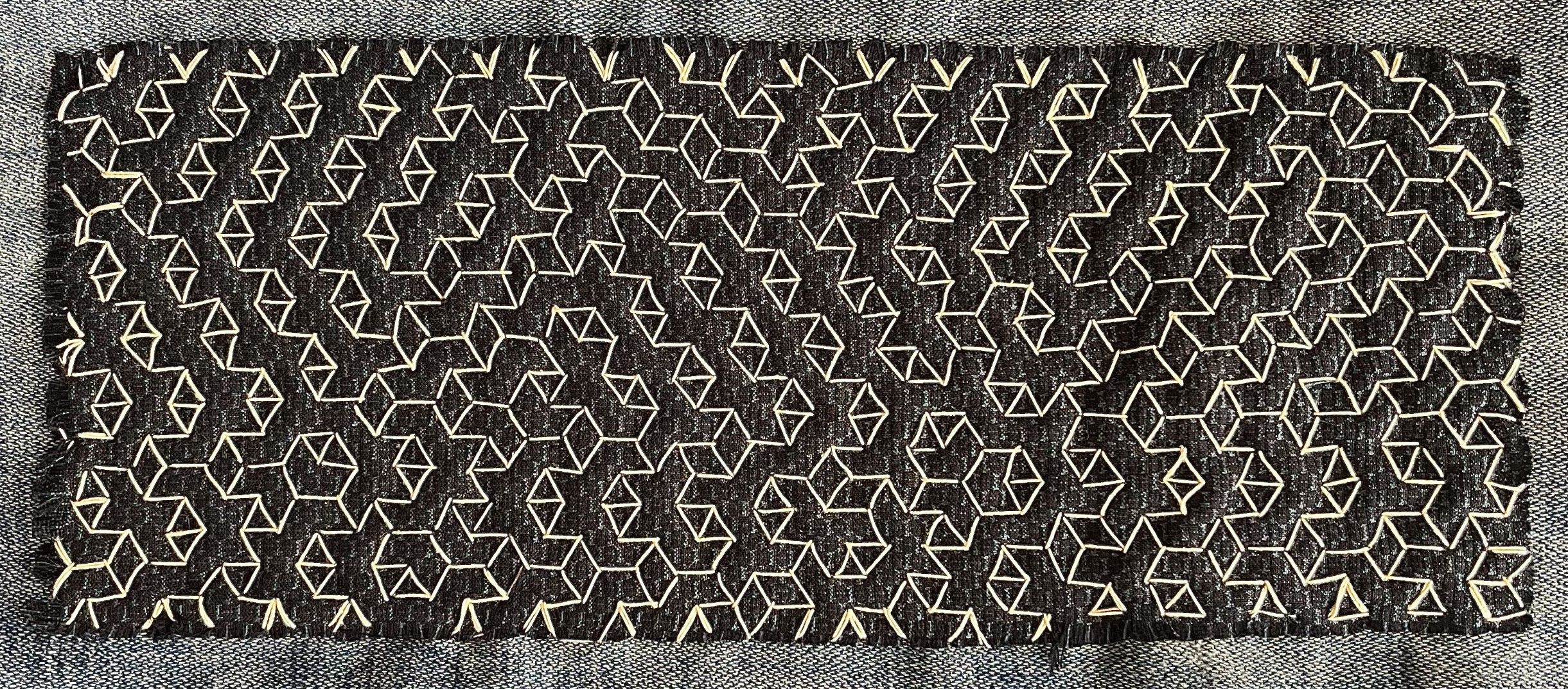}
	\caption{Fully packed hitomezashi on the isometric grid.\\ Stitched and photographed by Gwen Fisher and used with her kind permission.}
	\label{fig:1}
\end{figure}

For \textit{fully packed} hitomezashi on the square grid, that is, every line of running stitch in the two possible directions is present, each internal vertex has degree two when considered as part of the design. Previous discussions of hitomezashi on the equilateral triangular (isometric) grid have applied an analogous fully packed condition, thus making the degree of each vertex  \emph{three}  \cite{DKT}, \cite{Haran}. Figure \ref{fig:1} shows such stitching.

In this paper, I take a different approach. Rather than imposing the fully packed condition, only alternate lines of stitching in each direction of the three directions on the isometric grid are present. Each vertex visited by the stitching has degree two, and one quarter of the vertices are empty. To borrow a term from lattice models of statistical physics, let us call such stitching \textit{dilute}.

Kaplan previously adapted counted thread cross-stitch for the hexagonal/triangular grid \cite{Kaplan}. He used specially-sourced fabric with triaxial weave. When I first began to explore hitomezashi on the isometric grid, like the creator of Figure \ref{fig:1} I stitched on  patterned cotton that provided underlying hexagons \cite[p. 146]{Seaton}. My investigations recently received a huge push forward when I found a large piece of polyester mesh material (nearly 4 m$^2$) in a thrift shop, its holes forming the vertices of a triangular grid. Such knitted fabric is intended for sports garments and can be found with a variety of hole placements. It's slippery, and it lacks the precision of even weave linen or Aida cloth, but I found it much easier to work with for decorative hitomezashi than closely woven cotton. I have also made use of isometric grid paper to plan out designs.

\begin{figure}[h!tbp]
\centering
\begin{minipage}[b]{0.3\textwidth} 
	\includegraphics[width=\textwidth]{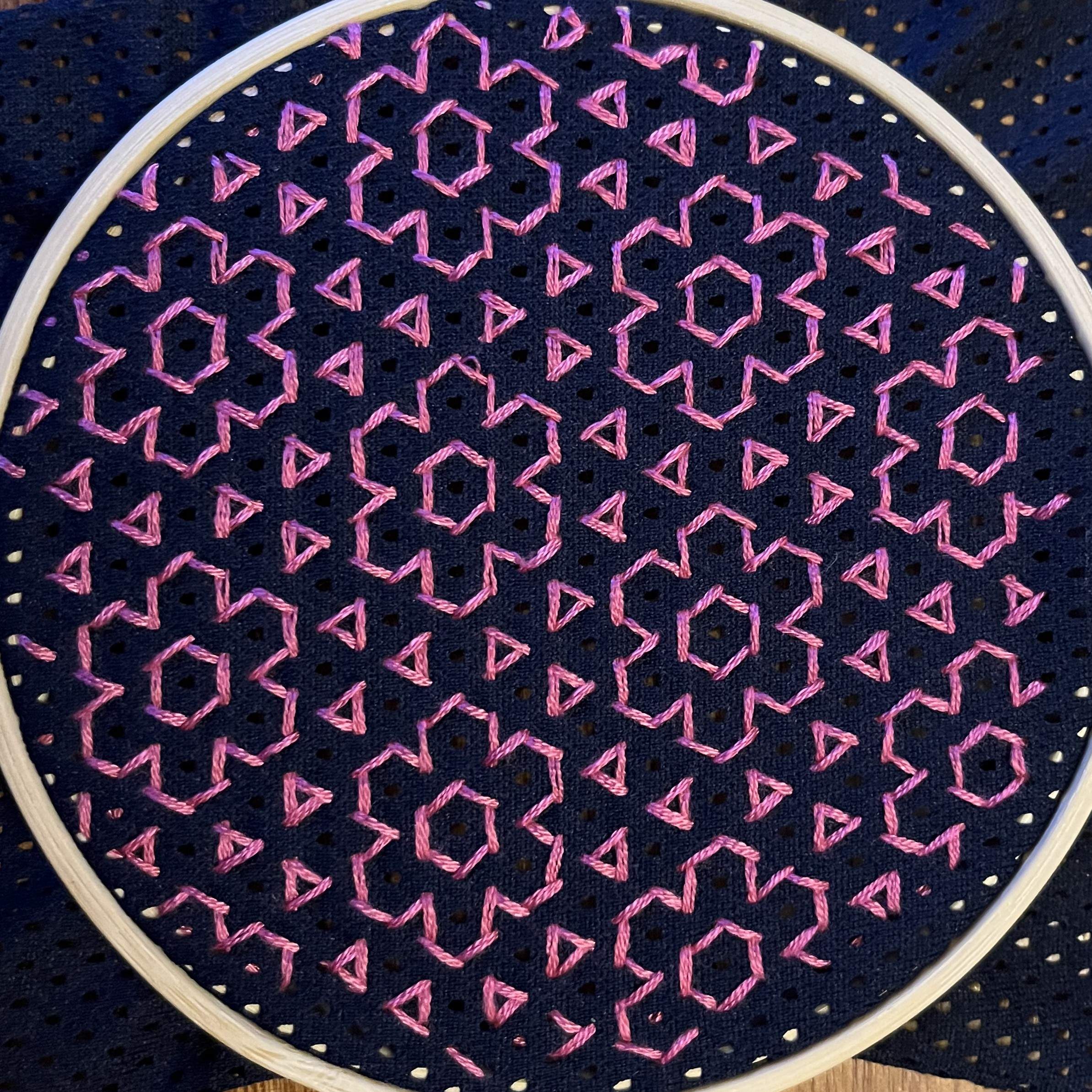}
        	\subcaption{} 
        	\label{fig:2a}
\end{minipage}
 ~	
\begin{minipage}[b]{0.3\textwidth} 
	\includegraphics[width=\textwidth]{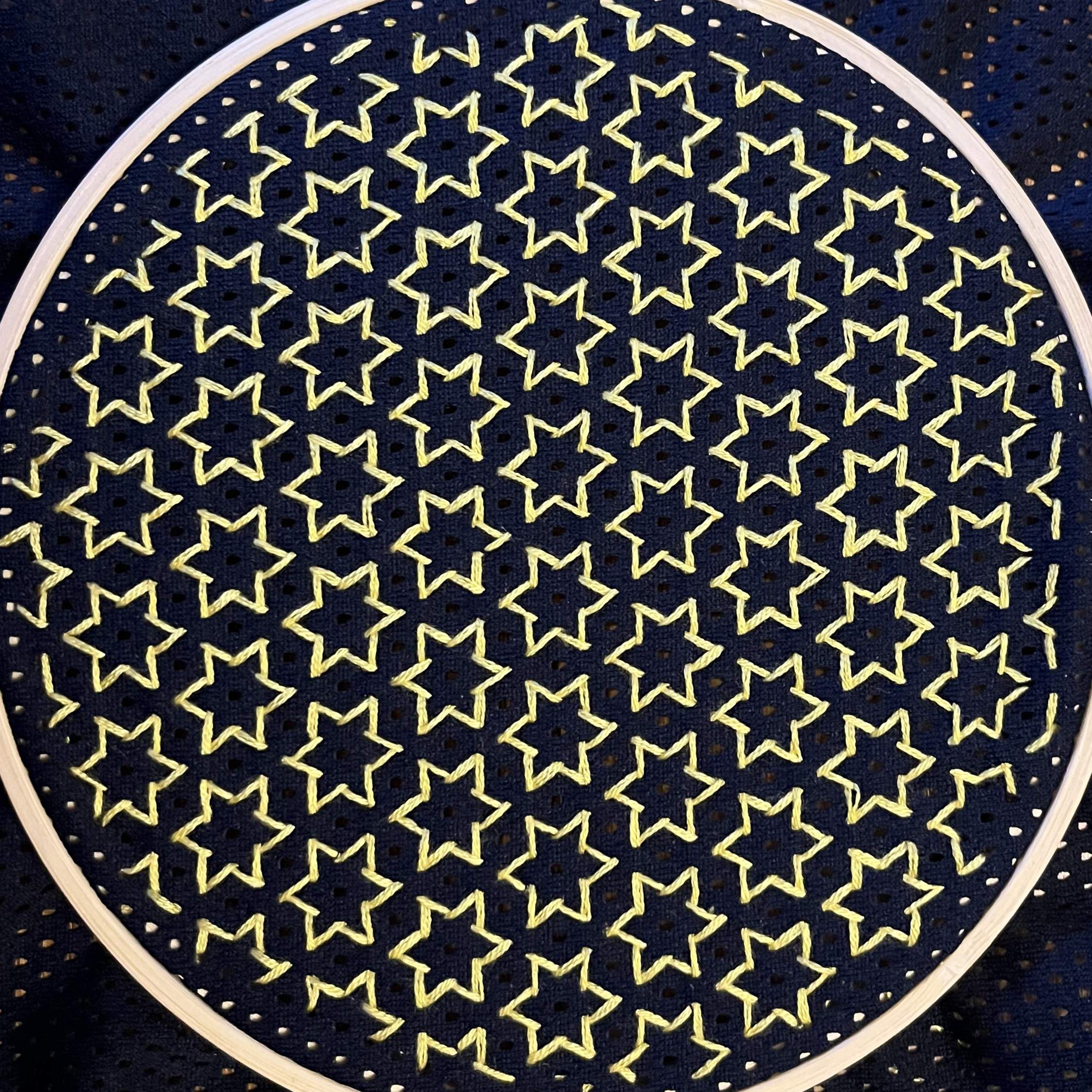}
        	\subcaption{} 
        	\label{fig:2b}
\end{minipage}
~
\begin{minipage}[b]{0.3\textwidth} 
	\includegraphics[width=\textwidth]{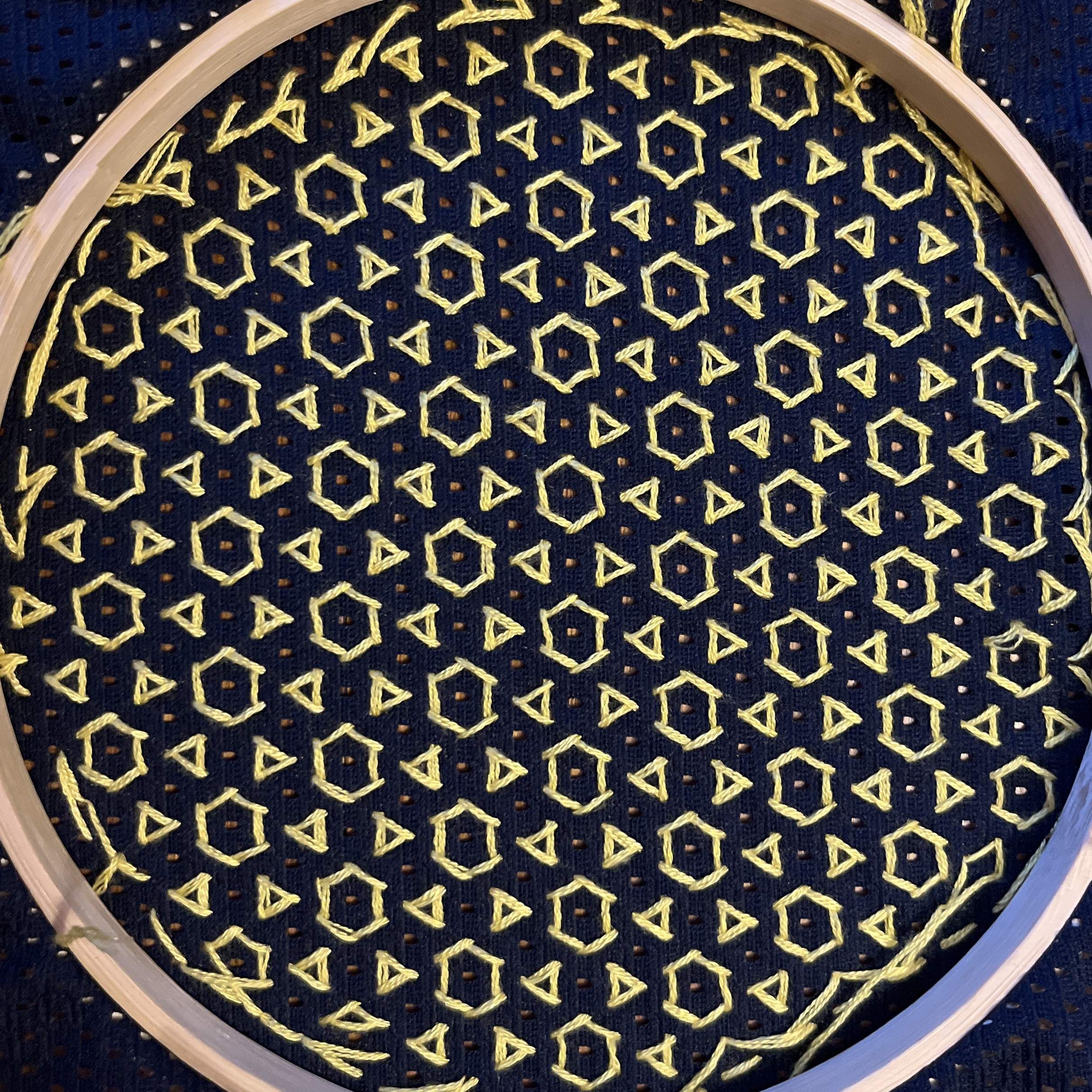}
        	\subcaption{} 
        	\label{fig:2c}
\end{minipage}

\caption{Dilute hitomezashi on the isometric grid.}
\label{fig:2}
\end{figure}

Unlike hitomezashi on the square grid, for which stitch dictionaries of myriad traditional patterns distill the work of countless domestic stitchers across centuries, with the isometric grid we find ourselves in uncharted waters albeit thus far only paddling. The one such design I have seen elsewhere, noticed on the Instagram account of the {\=O}tsuchi Sashiko Recovery Project, is similar to that of Figure \ref{fig:2a}. The six-petalled flowers were there worked as a frieze rather than a wallpaper and no name was given for the design. In the next section of this paper some designs with three- or six-fold symmetry are displayed. The final section demonstrates that the von Koch snowflake iterates, gloriously embellished with flowers and stars, can be stitched in hitomezashi on the isometric grid.

\section*{New Wallpaper Symmetries in Hitomezashi}

Martinez and Sen have proved which of the seventeen wallpaper patterns can be realised when the inherent constraints of fully packed square-lattice hitomezashi are applied \cite{Sen}. Necessarily, the five which feature angles of 60$^\circ$ or 120$^\circ$ do not appear.  Of the twelve patterns compatible with the square lattice, nine are possible in hitomezashi.  

In Figure \ref{fig:2b} what is arguably and perhaps surprisingly the simplest dilute design on the isometric grid is shown; all the lines of stitching in each direction are aligned. In a binary notation (see \cite{SeatHayes}) which indicates whether the first stitch in each line of stitching is present (1) or absent (0) on the front of the fabric, this design can be specified by the word $000\ldots$ in all three directions. The dual design to the repeated hexagrams is a pattern of triangles and hexagons shown in Figure \ref{fig:2c}. The appropriate wallpaper symmetry is p6mm for all three images in Figure \ref{fig:2}.

\begin{figure}[h!tbp]
\centering
\begin{minipage}[b]{0.3\textwidth} 
	\includegraphics[width=\textwidth]{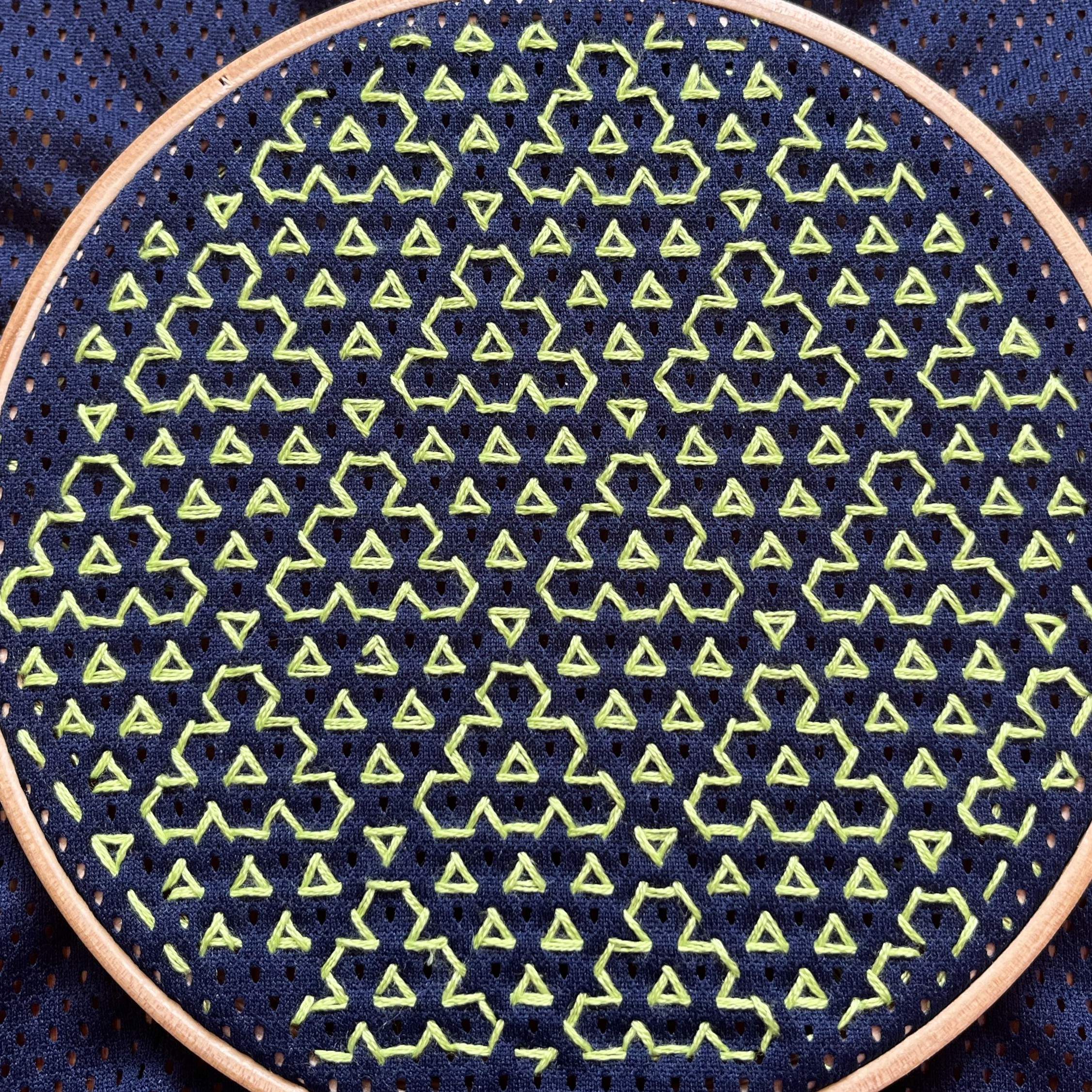}
        	\subcaption{} 
        	\label{fig:3a}
\end{minipage}
~
\begin{minipage}[b]{0.3\textwidth} 
	\includegraphics[width=\textwidth]{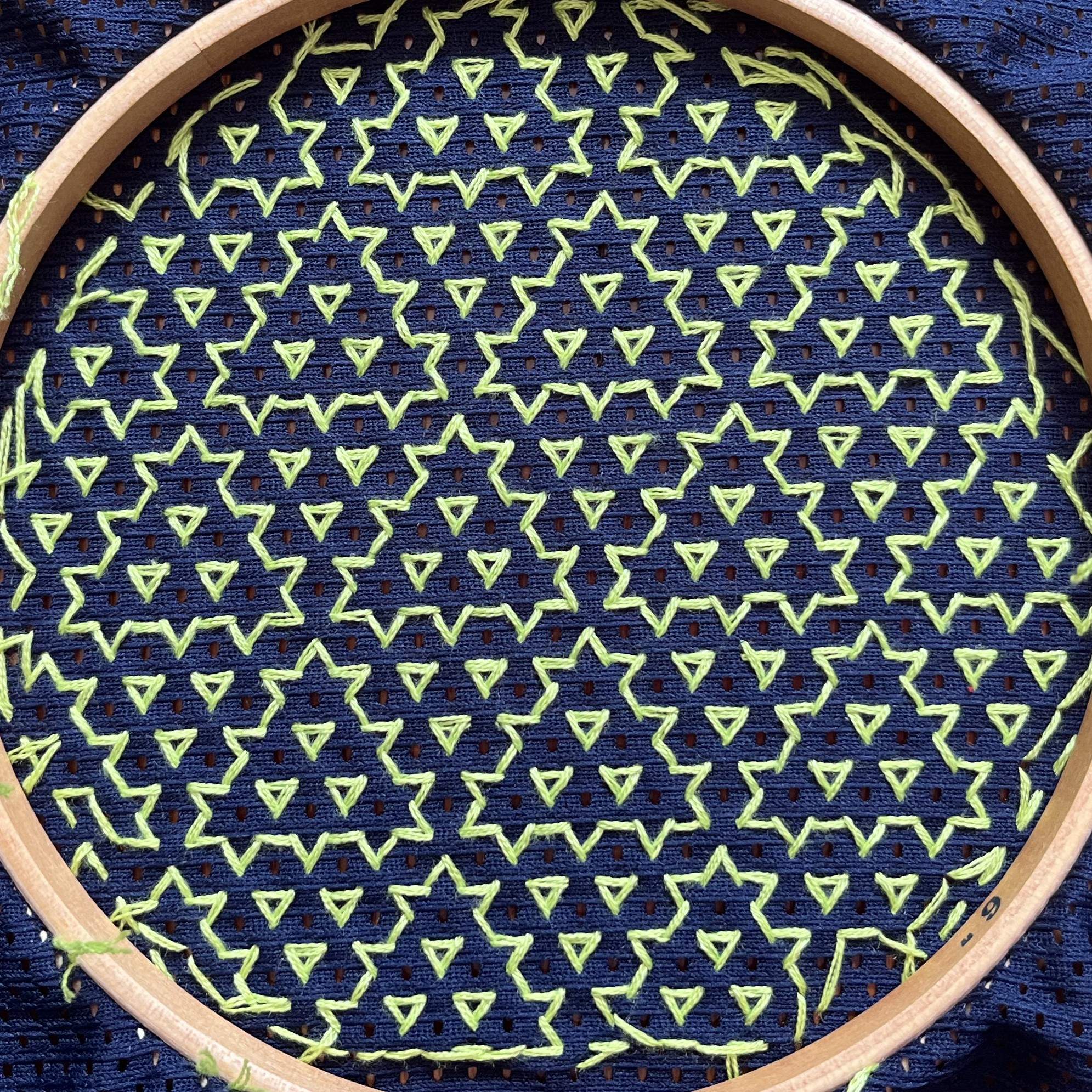}
        	\subcaption{} 
        	\label{fig:3b}
\end{minipage}
~
\begin{minipage}[b]{0.3\textwidth} 
\includegraphics[width=\textwidth]{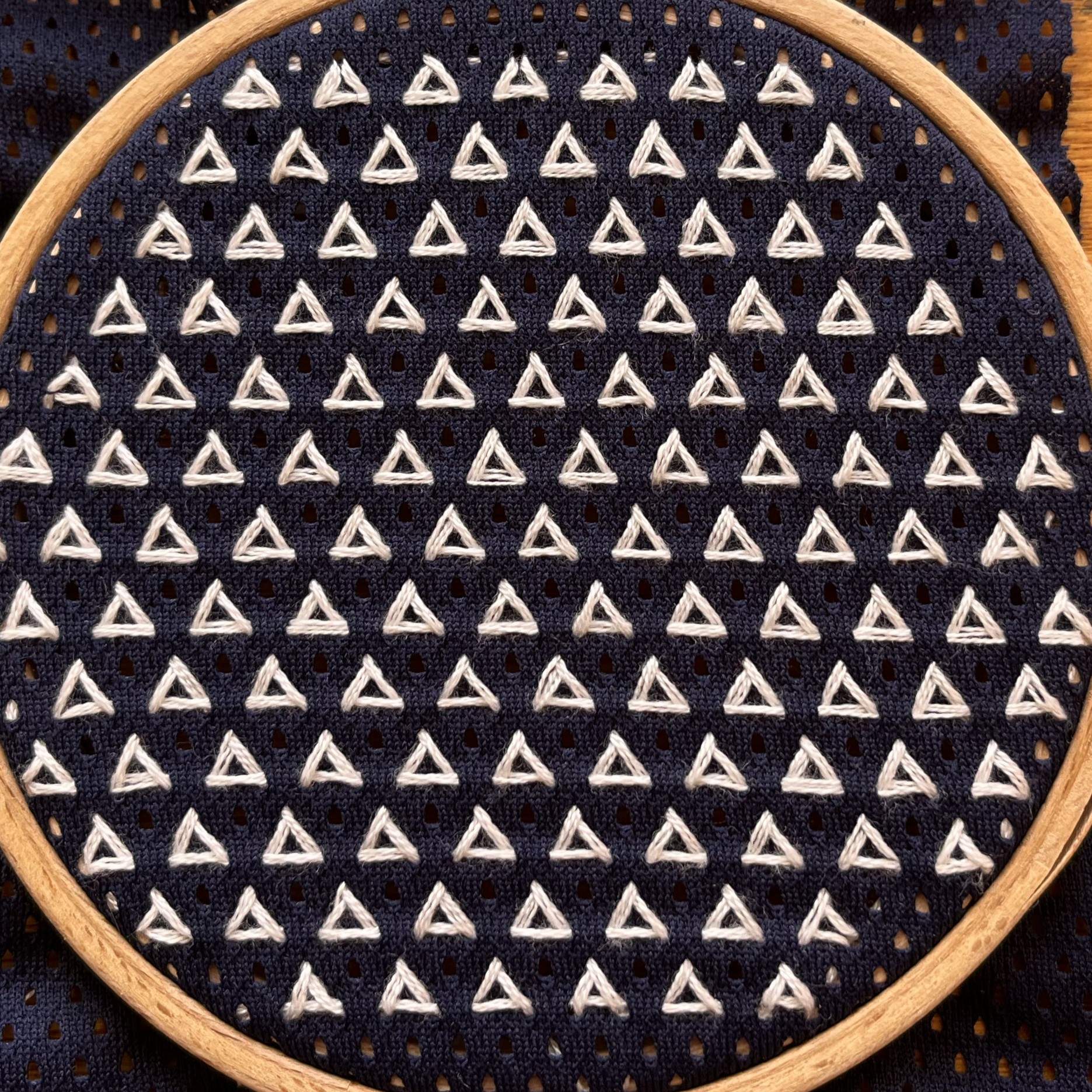}
        	\subcaption{} 
        	\label{fig:3c}
\end{minipage}

\caption{p3m1 patterns: (a), (b) a design and its dual; (c)  a self dual pattern. }
\label{fig:3}
\end{figure}

Figure \ref{fig:3} shows a design (and its dual) that manifests p3m1 symmetry. That is, there are centres of three-fold rotation as well as reflection axes passing through these centres. Three aligned lines of stitching are followed by one offset line, and then this sequence repeats: $00010001\ldots$. Again, the same stitching instructions are used in all three directions. Figure \ref{fig:3c} shows a design encoded as $0101\ldots$ which is self dual.

\section*{Koch Snowflake Iterates}

Previously, Carol Hayes and I discussed a recursive method to specify square lattice hitomezashi designs featuring the \textit{Fibonacci snowflakes} \cite{HayesSeat} and their generalizations \cite[Ch. 8]{Seaton}. The better known \textit{von Koch snowflake} fractal is the limit of a series of iterates generally considered to begin with an equilateral triangle, and then by applying an edge replacement rule, progressing through figures with six-fold symmetry and an increasingly crenelated boundary. I propose that the outlines of these iterates (apart from the triangle) arise as the result of specifying stitching patterns for dilute hitomezashi on the isometric lattice analogously.

\begin{figure}[h!tbp]
\centering
\begin{minipage}[b]{0.25\textwidth} 
	\includegraphics[width=\textwidth]{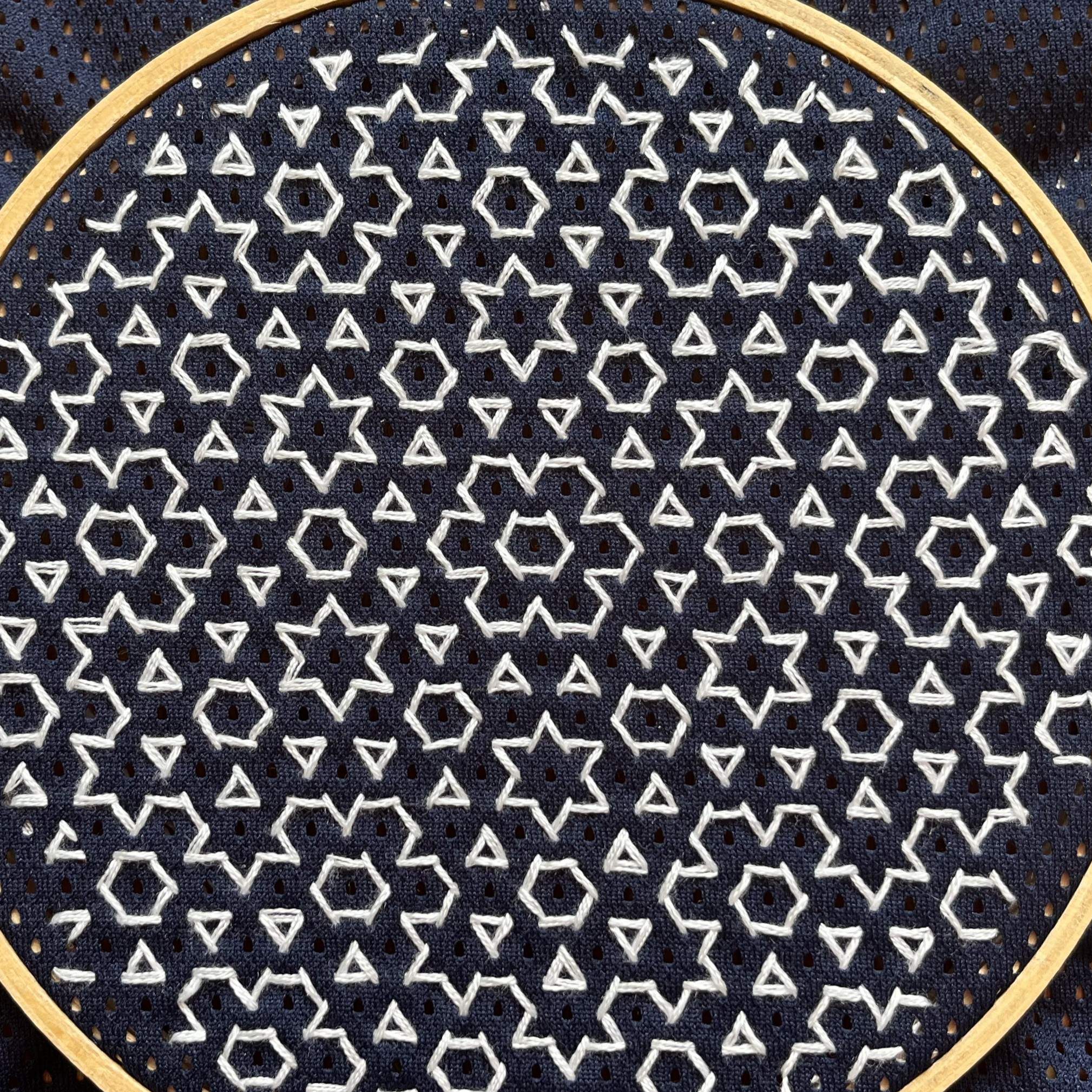}
        	\subcaption{} 
        	\label{fig:4a}
\end{minipage}
\qquad
\begin{minipage}[b]{0.25\textwidth} 
	\includegraphics[width=\textwidth]{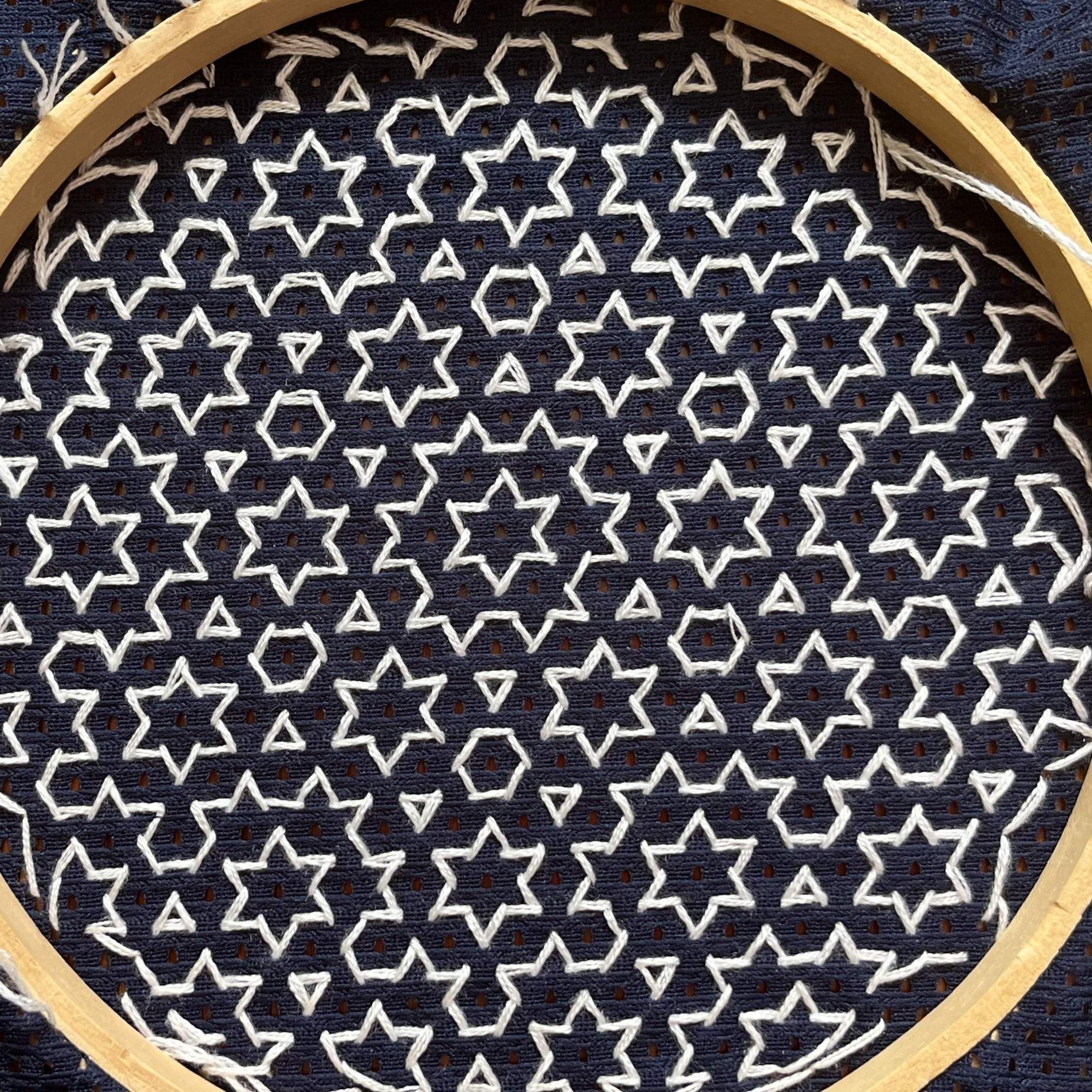}
        	\subcaption{} 
        	\label{fig:4b}
\end{minipage}

\caption{Koch snowflake iterate order $n=3$, and the dual design, filled with smaller motifs.}
\label{fig:4}
\end{figure}

\begin{figure}[h!tbp]
	\centering
	\includegraphics[width=5in]{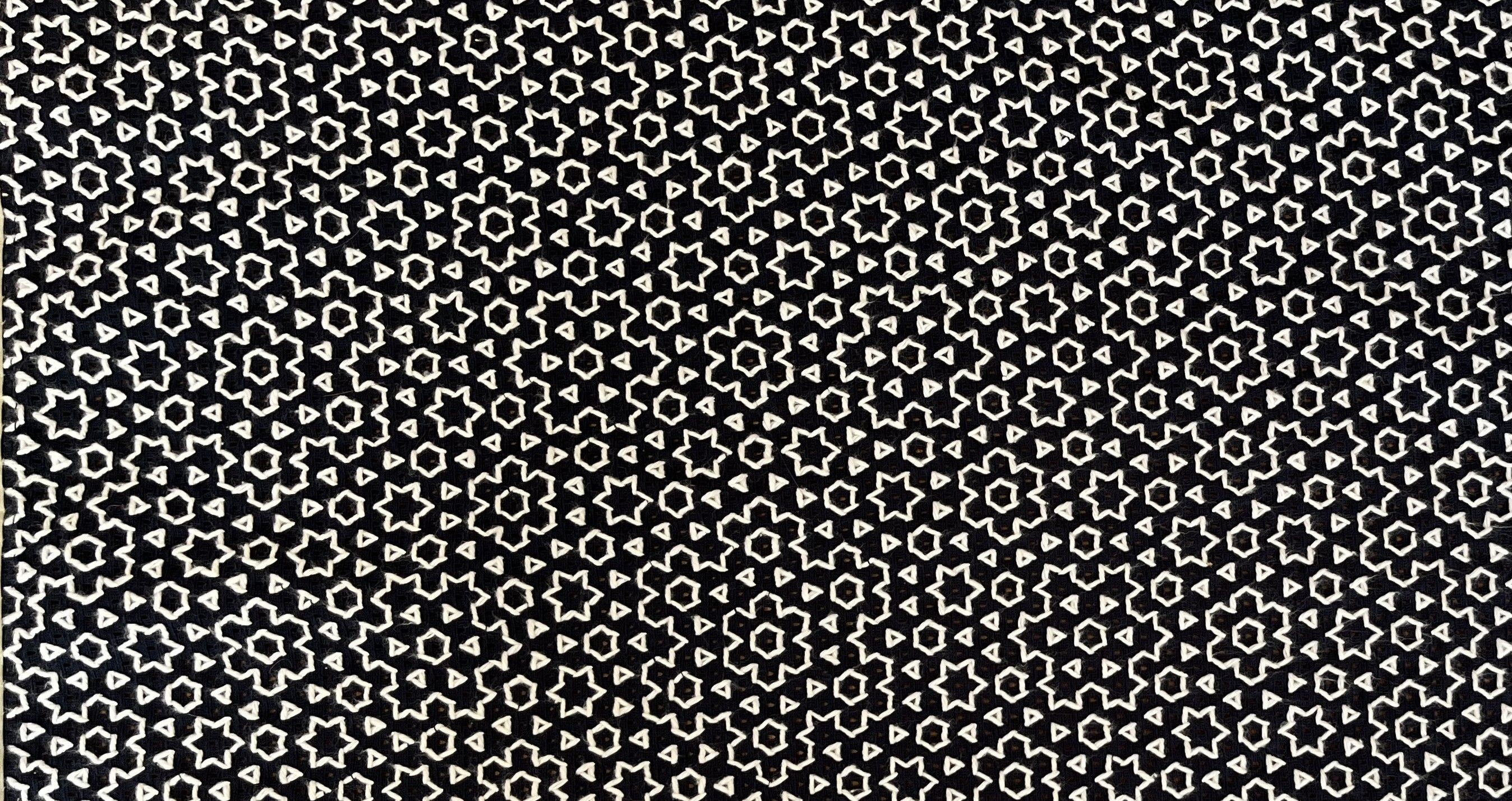}
	\caption{A detail from a larger piece \cite{Koch}, slightly more than half of an order 4 snowflake iterate.\\As well as other motifs, four order 2 snowflake iterates can be seen in this image.}
	\label{fig:5}
\end{figure}

Let $w$ be a binary word, and let $\overline{w}$ be the word obtained by the interchange $0\leftrightarrow 1$ in $w$, and $\widetilde{w}$ be obtained by reversing the order of the letters in $w$. For $n \in \mathbb N$, with $w_1=0$ let
\begin{equation}
 w_{n+1}= \widetilde{\overline{w_{n}}}\ \overline{w_n}\ w_n \label{RR}
\end{equation}
so that $w_2=110$ and $w_3=100001110$.  The result of using these words as stitching instructions in all three directions --- repeating them forwards and backwards $w_n\widetilde{w_n}$---  has been obtained up to order 4. The order 1 words (trivially) give the hexagram pattern of Figure \ref{fig:2b}. The higher order von Koch snowflake iterates do appear (see Figure \ref{fig:4a} and \ref{fig:5}), but this has the status only of an empirical observation and thus far I have no proof that Equation (\ref{RR}) `works' to all orders.  The words grow in length exponentially: $|w_n|=3^{n-1}$. Thus, to optimize the displayed detail, in Figure \ref{fig:5} just over half of the order 4 snowflake is shown.

In addition to the outlines of the Koch snowflakes, we see an array of structures within and between them. The recursive definition ensures that lower order iterates appear, but so too do six-petalled flowers and posies of such flowers. This seems to be a new way to look at---or inside---the von Koch snowflakes. Turning the hitomezashi over (see Figure \ref{fig:4b}), there are yet more motifs.

\section*{Summary and Conclusions}

This short paper serves only as an taster for  dilute hitomezashi on the isometric grid. While fully packed hitomezashi on this grid is hard to understand (see \cite{DKT}), the dilute version promises much. New wallpaper symmetries have been demonstrated, and fresh light thrown on the very familiar von Koch snowflakes.

    
{\setlength{\baselineskip}{13pt} 
\raggedright				

} 
   
\end{document}